\documentclass{svmult}

\usepackage{latexsym}
\usepackage{amssymb}
\usepackage{amsmath}
\usepackage{verbatim}
\usepackage{pdfsync}
\usepackage[curve,matrix,arrow,frame,tips]{xy}
\numberwithin{equation}{subsection}

\begin{document}

%%%%%%%%%%%%%%% This is the macros file for alll chapters of KYR
%%%%% last edited:  ssk  9/26/04

%%%%%%%%%%%%%%%  macros

%%%% mathbb

\newcommand{\A}{{\mathbb A}}
\newcommand{\C}{{\mathbb C}}
\newcommand{\F}{{\mathbb F}}
\newcommand{\G}{{\mathbb G}}
\newcommand{\R}{{\mathbb R}}
\newcommand{\Q}{{\mathbb Q}}
\newcommand{\X}{{\mathbb X}}
\newcommand{\Z}{{\mathbb Z}}
\newcommand{\HZ}{\widehat{\Z}}

%%%% mathrm

\newcommand{\rom}[1]{\text{\rm #1}}
\renewcommand{\roman}{\rm}

\newcommand{\Aut}{\text{\rm Aut}}
\newcommand{\CH}{\widehat{\text{\rm CH}}}
\newcommand{\cha}{{\text{\rm char}}}
\newcommand{\CHe}{\text{\rm CHeeg}}
\newcommand{\degh}{\widehat{\text{\rm deg}}}
\newcommand{\degH}{\widehat{\text{\rm deg}}}    %%% redundant def
\newcommand{\diag}{{\text{\rm diag}}}
\newcommand{\Diff}{\text{\rm Diff}}
\newcommand{\disc}{\text{\rm discr}}
\renewcommand{\div}{\text{\rm div}}
\newcommand{\divh}{\widehat{\text{\rm div}}}
\newcommand{\DS}{\text{\rm DS}}
\newcommand{\Ei}{\text{\rm Ei}}
\newcommand{\End}{\text{\rm End}}
\newcommand{\ev}{{\text{\rm ev}}}
\newcommand{\Gal}{\text{\rm Gal}}
\newcommand{\GL}{\text{\rm GL}}
\newcommand{\GSpin}{\text{\rm GSpin}}
\newcommand{\Hom}{\text{\rm Hom}}
\newcommand{\hor}{{\text{\rm horiz}}}
\newcommand{\id}{\text{\rm id}}
\newcommand{\im}{\text{\rm im}}
\renewcommand{\Im}{\text{\rm Im}}
\newcommand{\inv}{{\text{\rm inv}}}
\newcommand{\Jac}{\text{\rm Jac}}
\newcommand{\Leray}{{\mathrm L}}
\newcommand{\Lie}{\text{\rm Lie}}
\newcommand{\Mp}{\text{\rm Mp}}
\newcommand{\mult}{\text{\rm mult}}
\newcommand{\MW}{\text{\rm MW}}
\newcommand{\MWt}{\widetilde{\MW}}
\newcommand{\new}{\text{\rm new}}
\newcommand{\Nm}{\text{\rm Nm}}
\newcommand{\ord}{\text{\rm ord}}
\newcommand{\PGL}{\text{\rm PGL}}
\newcommand{\Pic}{\text{\rm Pic}}
\newcommand{\Pich}{\widehat{\text{\rm Pic}}}
\newcommand{\pr}{\text{\rm pr}}
\newcommand{\ra}{\text{\rm ra}}
\newcommand{\Rao}{\mathrm R}
\renewcommand{\Re}{\text{\rm Re}}
\newcommand{\sgn}{\text{\rm sgn}}
\newcommand{\sig}{\text{\rm sig}}
\newcommand{\SL}{\text{\rm SL}}
\newcommand{\SO}{\text{\rm SO}}
\newcommand{\Sp}{\text{\rm Sp}}
\newcommand{\Spec}{\text{\rm Spec}\, }
\newcommand{\Spf}{\text{\rm Spf}}
\newcommand{\supp}{\text{\rm supp}}
\newcommand{\Sym}{{\text{\rm Sym}}}
\newcommand{\tr}{\text{\rm tr}}
\newcommand{\type}{\text{\rm type}}
\newcommand{\Ver}{\text{\rm Vert}}
\newcommand{\vol}{\text{\rm vol}}
\newcommand{\Wald}{\text{\rm Wald}}

%%%% cals

\newcommand{\Cal}{\mathcal}     %%% this makes the old \Cal valid

\newcommand{\AHH}{\hat{\Cal A}}   % used??
\newcommand{\CHH}{\hat{\Cal C}}
\newcommand{\MM}{\Cal D}          % redefined!!!
\newcommand{\MMb}{\MM^\bullet}
\newcommand{\ssplit}{\text{\bf split}}
\newcommand{\whcc}{\widehat{\Cal C}}
\newcommand{\CO}{\mathcal O}
\newcommand{\COH}{\widehat{\CO}}
\newcommand{\M}{\Cal M}
\newcommand{\OB}{\Cal O_B}
\newcommand{\XX}{\mathcal X}
\newcommand{\bXX}{\bar\XX}
\newcommand{\wc}{\hat{\Cal C}}
\newcommand{\wch}{\wc^{\text{\rm hor}}}
\newcommand{\ZZ}{\Cal Z}
\newcommand{\ZH}{\widehat{\Cal Z}}   %%% redundant def's
\newcommand{\Zh}{\widehat{\Cal Z}}
\newcommand{\ZZh}{\ZZ^{\text{\rm hor}}}
\newcommand{\ZZv}{\ZZ^{\text{\rm ver}}}
\newcommand{\ZZhh}{\Zh^{\text{\rm hor}}}
\newcommand{\ZZhv}{\Zh^{\text{\rm ver}}}

%%%% math spacing

\newcommand{\nass}{\noalign{\smallskip}}
\newcommand{\snass}{\noalign{\vskip 2pt}}
\newcommand{\tent}[1]{ \vphantom{\vbox to #1pt{}} }   %%% !!!!

%%%% math fonts

\newcommand{\scr}{\scriptstyle}
\newcommand{\disp}{\displaystyle}

\font\cute=cmitt10 at 12pt
\font\smallcute=cmitt10 at 9pt
\newcommand{\kay}{{\text{\cute k}}}
\newcommand{\smallkay}{{\text{\smallcute k}}}

\renewcommand{\a}{\alpha}
\renewcommand{\b}{\beta}
\newcommand{\e}{\epsilon}
\renewcommand{\l}{\lambda}
\renewcommand{\L}{\Lambda}
\renewcommand{\o}{\omega}
\renewcommand{\O}{\Omega}
\renewcommand{\P}{\Phi}
\newcommand{\ph}{\varphi}
\newcommand{\phih}{\widehat{\phi}}
\newcommand{\wphi}{\widehat{\phi}}
\newcommand{\phit}{\widetilde{\phi}}
\newcommand{\s}{\sigma}
\newcommand{\vth}{\vartheta}

%%%% from Chapter VIII

%%\newcommand{\Gt}{\widetilde{G}}    %%%%  tilde's removed 6/20/04
%\newcommand{\Gt}{G}
%\newcommand{\Ph}{\Phi}
%\newcommand{\pht}{\widetilde{\phi}}
%%\newcommand{\Pht}{\widetilde{\Phi}}%%%%  tilde's removed 6/20/04
%\newcommand{\Pht}{\Phi}
%%\newcommand{\Pt}{\widetilde{P}}
%\newcommand{\Pt}{P}                 %%%%  tilde's removed 6/20/04
%%\newcommand{\Kt}{\widetilde{K}}    %%%%  tilde's removed 6/20/04
%\newcommand{\Kt}{K}
%%\newcommand{\It}{\widetilde{I}}    %%%%  tilde's removed 6/20/04
%\newcommand{\It}{I}
%\newcommand{\Jt}{\widetilde{J}}
%\newcommand{\lt}{\widetilde{\l}}
%\newcommand{\vp}{\varpi}
%

\newcommand{\Pt}{P}
\newcommand{\Ph}{\P}
\newcommand{\Pht}{\tilde \P}   %%%%%%%%    ****** conflicts *******
\newcommand{\Kt}{K}           %%%%  tilde's removed 6/20/04
\newcommand{\Mt}{M}
%%%%%%

%\newcommand{\Ph}{\Phi}      %%%%%%%%    ****** conflicts *******   temp %'ed
\newcommand{\pht}{\widetilde{\phi}}
\newcommand{\It}{I}
\newcommand{\Jt}{\widetilde{J}}
\newcommand{\lt}{\widetilde{\l}}
\newcommand{\vp}{\varpi}

\newcommand{\bom}{{\boldsymbol{\o}}}
\newcommand{\hbom}{\widehat{\bom}}
\newcommand{\ff}{{\bold f}}
\newcommand{\fsp}{\boldsymbol{f}_{\!\rm sp}}
\newcommand{\fev}{\boldsymbol{f}_{\!\rm ev}}
\newcommand{\fb}{\boldsymbol{f}}
\newcommand{\J}{\und{J}'}
\newcommand{\JJ}{\bold J'}
\newcommand{\V}{\bold V}
\newcommand{\xx}{\bold x}

\newcommand{\g}{{\mathfrak g}}
\renewcommand{\H}{\mathfrak H}

%%%%  math macros

\newcommand{\back}{\backslash}
\newcommand{\CT}[1]{\operatornamewithlimits{CT}_{#1}}
\renewcommand{\d}{\partial}
\newcommand{\db}{\bar\partial}
\newcommand{\dbar}{\bar{\partial}}
\newcommand{\gs}[2]{\langle \,#1,#2\,\rangle}
\newcommand{\Gt}{G}
\newcommand{\hfal}{h_{\text{\rm Fal}}}
\newcommand{\II}{\int^\bullet}
\newcommand{\isoarrow}{\ {\overset{\sim}{\longrightarrow}}\ }
\newcommand{\lisoarrow}{\ {\overset{\sim}{\longleftarrow}}\ }
\newcommand{\limdir}[1]{\underset{\underset{#1}{\rightarrow}}{\lim}}
\newcommand{\lan}{\operatorname{\langle}\hskip .5pt}
\newcommand{\ran}{\,\operatorname{\rangle}}
\newcommand{\lra}{\longrightarrow}
\newcommand{\doublelra}{\ {\overset{\scr\lra}{\scr\lra}}\ }
\newcommand{\nat}{\natural}
\newcommand{\notmid}{\mkern-5mu\not\mkern5mu\mid}
\newcommand{\Optoc}{\text{\rm Opt}(O_{c^2d},O_B)}
\newcommand{\psim}{\psi^{-}}
\newcommand{\qeq}{\ \overset{??}{=}\ }
\newcommand{\sh}{\sharp}
\newcommand{\thCH}{\theta^{\text{\rm ar}}}
\newcommand{\wht}{\widehat{\theta}}     %%% a replacement
\newcommand{\triv}{1\!\!1}
\renewcommand{\tt}{\otimes}
\newcommand{\und}[1]{\underline{#1}}
\newcommand{\z}{z}  %%% symbol used for the central sign 

\newcommand{\thMW}{\theta^{\text{\rm ar}}}
\newcommand{\tph}{\widetilde{\widehat\phi_1}}
\newcommand{\Pet}{\text{\rm Pet}}

%%%%% from Chapter IX

%%\newcommand{\Gt}{\widetilde{G}}    %%%%  tilde's removed 6/20/04
%\newcommand{\Gt}{G}
%\newcommand{\Ph}{\Phi}
%\newcommand{\pht}{\widetilde{\phi}}
%%\newcommand{\Pht}{\widetilde{\Phi}}%%%%  tilde's removed 6/20/04
%\newcommand{\Pht}{\Phi}
%%\newcommand{\Pt}{\widetilde{P}}
%\newcommand{\Pt}{P}                 %%%%  tilde's removed 6/20/04
%%\newcommand{\Kt}{\widetilde{K}}    %%%%  tilde's removed 6/20/04
%\newcommand{\Kt}{K}
%%\newcommand{\It}{\widetilde{I}}    %%%%  tilde's removed 6/20/04
%\newcommand{\It}{I}
%\newcommand{\Jt}{\widetilde{J}}
%\newcommand{\lt}{\widetilde{\l}}
%\newcommand{\vp}{\varpi}

%%%%%%%%%%%

%%%%  hacking

\newcommand{\thing}{ \raisebox{-6.4pt}{$\tilde{\tilde{}}$}  }   %%% some hacking from 4/14/04
\newcommand{\smallthing}{ \raisebox{-4.4pt}{$\scr\tilde{\tilde{}}$}  }
\newcommand{\ttilde}[1]{\overset{\smash{\thing}}{#1}}
\newcommand{\smallttilde}[1]{\overset{\smash{\smallthing}}{#1}}
\newcommand{\downhookarrow}{\hbox{$\downarrow\hskip -6.1pt\raisebox{6pt}{$\cap$}$}}

%%%% general formating

%\newcommand{\bysame}{\makebox[1.2cm][s]{\hrulefill ,\ }}   %%%%  \bysame is defined already
%\newcommand{\bysame}{$\underline{\text{\hbox to.5in{}}}$}   %%% was used in Textures
\providecommand{\bysame}{\makebox[3em]{\hrulefill}\thinspace}   %%% a fix:  cf. p 322 of Graetzer
\newcommand{\hfb}{\hfill\break}
\newcommand{\margincom}[1]{\marginpar{\bf\raggedright #1}}
\newcommand{\Sec}{\S}

%%%%%%%%%%%%%%

\numberwithin{equation}{section}
\setcounter{section}{0}
\setcounter{MaxMatrixCols}{15}

%%%%%%%%%%%%%%

\newtheorem{theo}{Theorem}[section]
\newtheorem{lem}[theo]{Lemma}
\newtheorem{prop}[theo]{Proposition}
\newtheorem{cor}[theo]{Corollary}
\newtheorem{conj}[theo]{Conjecture}
\newtheorem{rem}[theo]{Remark}      %%% seems not to exist in compositio.cls ???
\newtheorem{defn}[theo]{Definition}

\newcommand{\OO}{\text{\rm O}}
\newcommand{\UU}{\text{\rm U}}

\newcommand{\OK}{O_{\smallkay}}
\newcommand{\DI}{\mathcal D^{-1}}

\newcommand{\pre}{\text{\rm pre}}

\newcommand{\Bor}{\text{\rm Bor}}
\newcommand{\Rel}{\text{\rm Rel}}
\newcommand{\rel}{\text{\rm rel}}
\newcommand{\Res}{\text{\rm Res}}
\newcommand{\TG}{\widetilde{G}}

\newcommand{\PP}{\mathcal P}
\renewcommand{\OO}{\mathcal O}
\newcommand{\BB}{\mathbb B}
\newcommand{\GU}{\text{\rm GU}}
\newcommand{\Herm}{\text{\rm Herm}}

\newcommand{\FF}{\mathbb F}
\renewcommand{\MM}{\mathbb M}
\newcommand{\YY}{\mathbb Y}
\newcommand{\VV}{\mathbb V}
\newcommand{\LL}{\mathbb L}

\newcommand{\subover}[1]{\overset{#1}{\subset}}
\newcommand{\supover}[1]{\overset{#1}{\supset}}

\newcommand{\lcupover}[1]{{\scr #1}\,\cup\,\phantom{\scr #1}}
\newcommand{\rcupover}[1]{\phantom{\scr #1}\,\cup\,{\scr #1}}

\newcommand{\hgs}[2]{\{#1,#2\}}

\newcommand{\wh}[1]{\widehat{#1}}

\newcommand{\Ker}{\text{\rm Ker}}
\newcommand{\YYbar}{\overline{\YY}}
\newcommand{\MMbar}{\overline{\MM}}
\newcommand{\oY}{\overline{Y}}
\newcommand{\dra}{\dashrightarrow}
\newcommand{\dlra}{\longdashrightarrow}

\newcommand{\yy}{\text{\bf y}}

\newcommand{\red}{\text{\rm red}}

\newcommand{\inc}{\text{\rm inc}}

\newcommand{\OKs}{O_{\smallkay,s}}
\newcommand{\OKr}{O_{\smallkay,r}}
\newcommand{\Xs}{X^{(s)}}
\newcommand{\Xo}{X^{(0)}}

\newcommand{\Nilp}{\text{\rm Nilp}}
\newcommand{\beq}{\begin{equation}}
\newcommand{\eeq}{\end{equation}}

\newcommand{\uA}{\und{A}}
\newcommand{\cutter}{\medskip\medskip \hrule \medskip\medskip}

\newcommand{\YYb}{\bar{\YY}}
\newcommand{\Wk}{W_{\smallkay}}
\renewcommand{\ss}{\text{\rm ss}}
\newcommand{\Jor}{\text{\rm Jor}}
\newcommand{\WK}{W_{\smallkay}}

\newcommand{\da}{{[r]}}
\newcommand{\darr}{\dagger}
\newcommand{\sha}{\sharp}   

\newcommand{\oo}[1]{\overset{\ o}{#1}}
\newcommand{\ua}{\underline{A}}
\newcommand{\um}{\underline{M}}
\newcommand{\ud}{\underline{d}}

\newcommand{\Y}{\mathbb Y}
\newcommand{\Yb}{{\bar\Y}}
\newcommand{\xibold}{\text{\boldmath$\xi$\unboldmath}}

\newcommand{\MY}{M(\bar{\mathbb Y}^r)}

\renewcommand{\top}{\text{\rm top}}

\newcommand{\OC}{O_C}

\title*{A note about special cycles on moduli spaces of K3 surfaces}
\author{Stephen Kudla}
\institute{Stephen Kudla \  Department of Mathematics, University of Toronto, 40 St George St. BA6290, 
Toronto, ON M5S 2E4, Canada \email{skudla@math.toronto.edu}}

\maketitle

\date{today}

\centerline{\bf Abstract\hfill}
\medskip

We describe the application of the results of Kudla-Millson on the modularity of 
generating series for cohomology classes of special cycles to the case of lattice polarized K3 surfaces.
In this case, the special cycles can be interpreted as higher Noether-Lefschetz  loci. 
These generating series can be paired with the cohomology classes of complete subvarieties of the moduli 
space to give classical Siegel modular forms with higher Noether-Lefschetz numbers as Fourier 
coefficients.  Examples of such complete families associated to quadratic spaces over totally real number fields are constructed.

\vskip .5in

\centerline{\bf Introduction\hfill}

\medskip

This article contains a short survey of some results about special cycles on certain Shimura varieties that 
occur as moduli spaces of lattice polarized K3 surfaces.   
The two points that may be
of interest are the Siegel modular forms arising as generating series for higher Noether-Lefschetz numbers
and the description of some complete subvarieties in certain of these moduli spaces. 
The families of K3 surfaces parametrized by these subvarieties ought to have particularly nice properties and I do not know to what extent they already 
occur implicitly or explicitly in the literature.
Finally, I ask the indulgence of the experts in this area for my very naive treatment of things that may be very well known to them.

\section{Special cycles for orthogonal groups}

We begin by reviewing a very special case of old joint work with John Millson, \cite{kmI}, \cite{kmII}, \cite{KM.IHES}.

\subsection{} Let $L$, $(\ ,\ )$ be a lattice with a $\Z$-valued symmetric bilinear form of signature $(2,n)$. In particular,
for   $V=L_\Q=L\tt_\Z\Q$, the dual lattice
$$L^\vee = \{\  x\in V \mid (x,L) \subset \Z\ \}$$ 
contains $L$.  Let 
 \begin{align*}
D=D(L)&= \{\ w\in V_\C\mid (w,w)=0, \ (w, \bar w)>0\, \}/\C^\times\\
\nass
{}&\simeq \{ \text{oriented positive $2$-planes in $V_\R$}\}\\
\nass
{}&\simeq \SO(V_\R)/K
\end{align*}
be the associated symmetric space, where $K\simeq SO(2)\times SO(n)$ is the stabilizer of 
an oriented positive $2$-plane.  Here $V_\C =V\tt_\Q\C$ (resp. $V_\R = V\tt_\Q \R$), and we sometimes write
$\und{w}$ for the image of $w\in V_\C$ in $D$. 
Let $\Gamma_L\subset \SO(V)$ be the isometry group of $L$
and let 
$\Gamma\subset \Gamma_L$ 
be a subgroup of finite index. 
Then 
$$M_{\Gamma}= \Gamma\back D(L),\qquad \dim_\C M_\Gamma = n$$ 
is (isomorphic to) a quasi-projective variety.  This variety is a connected component of a Shimura variety\footnote{For $n\le 2$, we need to add the 
condition that $\Gamma$ be a congruence subgroup; this is automatic for $n\ge 3$. In general, $M_\Gamma$ can have $2$ components, since $D(L)$ does. }
and has a model defined over a cyclotomic field. 
For small values of $n$, this space can be the moduli space of polarized K3 surfaces and, for 
smaller values, of abelian varieties.

\subsection{}  To define special cycles, suppose that $x\in L$ is a vector with $(x,x)<0$, and let
$$D_x= \{ \und{w}\in D\mid (x,\und{w})=0\} = D(L\cap x^\perp).$$
Thus $D_x$ has codimension $1$ in $D$ and gives rise to a divisor 
$$Z(x): \Gamma_x\back D_x\lra \Gamma\back D = M_\Gamma,$$
in $M_\Gamma$, where $\Gamma_x$ is the stabilizer of $x$ in $\Gamma$. 
We call such a divisor, which  depends only on the $\Gamma$-orbit of $x$, a special divisor. 

Composite divisors are defined as follows.   Let 
$\Gamma_L^o\subset \Gamma_L$ be the subgroup of isometries that 
act trivially on $L^\vee/L$, and suppose that $\Gamma\subset \Gamma_L^o$. 
For $t\in \Z_{>0}$ and $h\in L^\vee$ a coset representative for $L^\vee/L$,
let 
$$L(t,h)= \{ x\in L+h \mid (x,x) = -t\}.$$
There is an associated special divisor
$$Z(t,h)= \sum_{\substack{x\in L(t,h)\\ \snass \mod \Gamma}}  Z(x).$$

\noindent{\bf Remark.} For the following observation, see 
\cite{vG}, Lemma 1.7.  The space $D$ parametrizes polarized Hodge structures of weight $2$ on the rational vector space $V$ with $\dim_\C V_\C^{2,0}=1$. 
Such a HS is {\it simple} if it does not contain any proper rational Hodge substructure. Then, in fact,  the polarized HS corresponding to $\und{w}\in D$ is simple
if and only if $\und{w} \notin D_x$ for any nonzero vector $x\in V$. Thus 
the set  
$$D\ - \bigcup_{x\in V, \ x\ne 0} D_x$$
parametrizes the simple HS's of this type. 

\medskip

More generally, for $1\le r \le n$, consider an $r$-tuple of vectors
\begin{align*}
\xx&=[x_1,\dots, x_r], \qquad x_i\in V,\\
\noalign{\noindent and suppose that}
T(\xx)&= -(\xx,\xx) = -((x_i,x_j)) >0.\\
\noalign{\noindent Let }
D_{\xx}&=  \{ \und{w}\in D\mid (\xx,\und{w})=0\},\\
\nass
\Gamma_{\xx}&=\text{stablizer of $\xx$ in $\Gamma$},\\
\noalign{\noindent and}
Z(\xx):& \ \ \Gamma_{\xx}\back D_{\xx}\lra \Gamma\back D = M_\Gamma.
\end{align*}
The condition $T(\xx)>0$ implies that $D(\xx)$ has codimension $r$ in $D$ so that the 
special cycle $Z(\xx)$ has codimension $r$ in $M_\Gamma$.  Again, this cycle depends only on the $\Gamma$-orbit of $\xx$.

Again there is a composite version.  For 
\begin{align*}
T&\in \Sym_r(\Z)_{>0}, \qquad \bold h\in (L^\vee)^r,\\
\noalign{\noindent let}
L(T, \bold h)&= \{ \xx\in L^r +\bold h\mid T(\xx)=T\},\\
\noalign{\noindent and define the cycle}
%\nass
Z(T,\bold h)&= \sum_{\substack{\xx\in L(T,\bold h)\\ \snass\mod \Gamma}}  Z(\xx).
\end{align*}
These are the special cycles in question. 

\medskip

\noindent{\bf Remarks.}
a) 
In the rest of this note, we will suppress the coset parameter $h$, although it plays an evident and important role in 
many places in the literature. \hfb
b) We can make the same construction with any $T\in \Sym_r(\Z)_{\ge0}$ except that we require that 
for $\xx\in L(T,\bold h)$ the subspace spanned by the components of $\xx$ has dimension equal to the rank of $T$. 
We write $Z^{\text{naive}}(T,\bold h)$ for the resulting cycle. It has codimension equal to the rank of $T$. For example, for $\bold h=0$,  
$$Z^{\text{naive}}(0,\bold h) = M_\Gamma,$$
while, for $\bold h\ne 0$,  $Z^{\text{naive}}(0,\bold h)$ is empty. 
\medskip

\section{Modular generating series}

First, for  $\tau=u+iv \in \H$, the upper half-plane, consider the series
$$\phi_1(\tau) = [Z(0)]+\sum_{t\in\Z_{>0}}  [Z(t)]\, q^t, \qquad [Z(t)] \in H^2(M_\Gamma).$$
Here $[Z(0)] = [\omega] = c_1(\Cal L)$ is the Chern class of the tautological line bundle 
$$\Cal L \lra  M_\Gamma,$$ 
defined by: 
$$\Cal L= b^*(\Cal O(-1)), \qquad b: D(L) \lra D(L)^\vee,$$
where 
$$D^\vee(L) =\{\ w\in V_\C\mid (w,w)=0\ \}/\C^\times  \subset \mathbb P(V_\C)$$
is the compact dual of $D(L)$. Here $H^r(M_\Gamma)$ is the usual Betti cohomology group of the (quasi-projective) variety $M_\Gamma$ with complex 
coefficients\footnote{If $\Gamma$ has fixed points on $D(L)$, $M_\Gamma$ is viewed as an orbifold and $H^r(M_\Gamma)$ 
is the space of $\Gamma/\Gamma_1$-invariants in $H^r(M_{\Gamma_1})$ where $\Gamma_1\subset \Gamma$ is a normal subgroup 
of finite index which acts freely on $D(L)$.}. 

\begin{theo}[\cite{KM.IHES}]
$\phi_1(\tau)$ is an elliptic modular form of weight $\frac{n}2+1$ and level determined\footnote{
This means that the components of $\phi_1(\tau)$ with respect to any basis of the finite dimensional space $H^2(M_\Gamma)$ 
are scalar valued modular forms of the given weight and level. The level divides $4 |L^\vee/L|$ and  
is determined by the usual recipe for theta functions, cf., for example, \cite{shimura}, section 2.} by $L$, 
valued in $H^2(M_\Gamma)$. 
\end{theo}

\medskip
The analogous generating series with values in the first Chow group $\text{\rm CH}^1(M_\Gamma)\tt_\Z\C$ of $M_\Gamma$ was considered by Borcherds\footnote{
More precisely, he views $M_\Gamma$ as an orbifold/stack and defines a group by generators and relations that maps to the usual Chow group, at least 
after tensoring with $\C$.}

\begin{theo}[\cite{Bo2}] 
$$\phi_{1}^{CH}(\tau) = \{Z(0)\}+\sum_{t\in\Z_{>0}}  \{Z(t)\}\, q^t, \qquad \{Z(t)\}\in \text{\rm CH}^1(M_\Gamma)\tt\C$$
is an elliptic modular form of weight $\frac{n}2+1$, etc. Its image under the cycle class map
$$\text{\rm CH}^1(M_\Gamma)\tt\C \lra H^2(M_\Gamma)$$
is $\phi_1(\tau)$.  
\end{theo}

More generally, for any $r$, $1\le r\le n$, and for $\tau = u+iv \in \H_r$,  the Siegel space of genus $r$, 
we can define a generating series
\begin{align*}
\phi_r(\tau) &= \sum_{  \substack{T\in \Sym_r(\Z)\\ \snass T\ge 0}  } [Z^{\text{naive}}(T)]\cup [\o]^{r-\text{rank} T}\, q^T,\\
\nass
{}&= [\o]^r + \sum_{\text{rank} T<r}   [Z^{\text{naive}}(T)]\cup [\o]^{r-\text{rank} T}\, q^T + \sum_{T\in \Sym_r(\Z)_{> 0}} [Z(T)]\, q^T.
\end{align*}
The point is that we have shifted the classes $[Z^{\text{naive}}(T)]$ by a suitable power of $[\o]$ so that all of the coefficients lie in $H^{2r}(M_\Gamma)$. 

\begin{theo} [\cite{KM.IHES}] {\rm (contd.)}
$\phi_r(\tau)$ is a Siegel modular form of weight $\frac{n}2+1$ and level determined by $L$, 
valued in $H^{2r}(M_\Gamma)$. 
\end{theo}

\begin{rem}
{\bf(1)}
In \cite{kudla.duke} I asked whether the Chow group version $\phi_{r}^{CH}(\tau)$ is a Siegel modular form valued in $\text{\rm CH}^r(M_\Gamma)$. 
 Using Borcherds' result mentioned above and an inductive argument based on Fourier-Jacobi expansions, 
Wei Zhang proved this in his Columbia thesis, \cite{wei.zhang}, conditionally on some finiteness/convergence assumption. \hfb
{\bf(2)} A key point is that the pairing of $\phi_r(\tau)$ with any compactly supported class $c\in H^{2n-2r}_c(M_\Gamma)$ defines a classical Siegel modular form\footnote{This form 
will be vector-valued if we keep track of the parameter $\bold h\in L^\vee/L$.} of the same weight and level as $\phi_r(\tau)$. \hfb
{\bf(3)} It is worth noting  that the modularity of the $\phi_r(\tau)$ is proved by constructing a theta function 
$\theta(\tau,\ph_{KM}^{(r)})$ valued in $A^{(r,r)}(M_\Gamma)$ the space of smooth differential forms of type $(r,r)$ on $M_\Gamma$. 
This theta function is a non-holomorphic modular form of weight $\frac{n}2+1$ in $\tau$, analogous to the classical Siegel theta function for an 
indefinite quadratic form. Moreover, it is a closed $(r,r)$-form and its cohomology class $[\theta(\tau,\ph_{KM}^{(r)})]\in H^{(r,r)}(M_\Gamma)$ 
is $\phi_r(\tau)$.  In particular, if the class $c\in H^{2n-2r}_c(M_\Gamma)$ is the class of an algebraic cycle $c = [S]$, 
as will be the case in several examples below, the pairing of $c$ and $\phi_r(\tau)$ is given by an integral
$$c \cdot \phi_r(\tau) = \int_{S} \theta(\tau,\ph_{KM}^{(r)}).$$
{\bf(4)} The extension of such integrals to classes $c$ that are not compactly supported has been studied in important work of 
Funke and Millson, \cite{funke.millson.I}, \cite{funke.millson.II}, \cite{funke.millson.III}, and involves interesting correction terms related to the boundary of 
the Borel-Serre compactification of $M_\Gamma$. 
\end{rem}

\section{The case of K3 surfaces} 
In some cases, the $M_\Gamma$'s and the $Z(T)$'s can be interpreted in terms of moduli spaces of lattice polarized\footnote{A good references for lattice polarized K3 surfaces are \cite{dolgachev.mirror} and  \cite{dolgachev.kondo.istanbul}. } K3 surfaces. Here is an amateur's version of this.
Start from the K3 lattice:
$$K= H^3\oplus E_8(-1)^2, \qquad \sig(K)=(3,19),$$
where $H$ is the unimodular hyperbolic plane, and, letting $K_\Q = K\tt_\Z\Q$, choose an orthogonal decomposition
$$K_\Q= V\oplus V'$$
with 
$$\sig(V')= (1,19-n),\qquad \sig(V) = (2,n),$$
for some $n$, $0\le n\le 19$. 
Let 
$$L= K\cap V, \qquad L'= K\cap V',$$
be the corresponding primitive sublattices in $K$.  Note that they are both even integral. %and that

A marked $L'$-polarized K3 surface is a collection $(X,u,\l)$ where $X$ is an algebraic K3 surface over $\C$,
$$u: H^2(X,\Z)\isoarrow K$$
is an isometry for the intersection form on $H^2(X,\Z)$,  a marking, and 
$$
\l: L'\hookrightarrow  \Pic(X)\subset H^2(X,\Z)
$$
is an embedding such that
$$
u\circ \l: L'\hookrightarrow K
$$
is the given inclusion.
Moreover, $\l$ is required to satisfy the `ample cone' condition\footnote{More precisely, 
let $(L')_{-2}$ be the set of vectors $x$ in $L'$ with $(x,x)=-2$, and let $\Cal V(L')^0$ be one component of 
the set of vectors $x\in L'_\R$ with $(x,x)>0$, i.e., a choice of positive light cone. Finally, let 
$C(L')^+$ be a component of 
$$\Cal V(L')^0 - \bigcup_{x\in (L')_{-2}} x^\perp.$$
Then the ample cone condition is that 
$$\l_\R(C(L')^+) \cap \Cal K_X \ne \emptyset$$
where $\Cal K_X$ is the closure of the ample cone in $\Pic(X)_\R$. }, cf. section 10 of \cite{dolgachev.kondo.istanbul}. 

These conditions imply that 
\begin{align*}
(i)\qquad &\text{the embedding $\l:L'\hookrightarrow \Pic(X)$ is primitive and isometric}\\
\nass
(ii)\qquad&\text{the period point of $(X,u)$ lies in $D(V)$.}
\end{align*} 
Recall that the period point of $(X,u)$ is the complex line $u_\C(H^{2,0}(X))$, where \hfb
$u_\C: H^2(X,\Z)\tt_\Z\C \lra K_\C$ is the complex linear extension of $u$. 
The moduli space of such gadgets $(X,\l)$, obtained by eliminating the marking, is the quotient
$$M_\Gamma = \Gamma\back D(V),$$
$$\Gamma = \{\gamma \in \Gamma_L \mid \gamma\vert_{ L^\vee/L} = 1 \}.$$
This gives a moduli interpretation of $M_\Gamma$.  Details of this construction can be found in \cite{dolgachev.kondo.istanbul}, section 10. 

\begin{rem}\label{rem3.1} {\bf(1)} It is easy to check that, for $n\le 17$, any {\it rational} quadratic space of signature $(2,n)$ can occur as $V$.
If $n=18$ or $19$, there are restrictions on $V$, since then $L'_\Q$ has rank $2$ or $1$.\hfb %, see below.\hfb
{\bf(2)} 
In certain cases, a precise description of what lattices $L'$'s and $L$'s can occur was given by Nikulin \cite{nikulin.quad}.  As summarized in
\cite{morrison}, section 2, the result is the following. If $0\le n\le 9$, then any even integral lattice $L$ of signature $(2,n)$ can occur, 
and, if $n<9$, the primitive embedding $L\hookrightarrow K$ is unique up to an isometry of $K$. 
Similarly, if $0\le n'\le 10$, then any even integral lattice $L'$ of signature $(1,n')$ can occur, and, if $n'<10$, the embedding $L'\hookrightarrow K$, 
is unique up to an isometry of $K$.   
\hfb
{\bf(3)} In such a family the generic element has Picard number $\rho(X) = \text{rank}(L' )=20-n$, 
\end{rem}

\medskip

For small values of $n$, the resulting $M_\Gamma$'s are familiar classical objects. Here is a little table:

\begin{center}
\begin{tabular}{| l | c | c | l | l |}
\hline
$n$ & $\rho$ & $G=\SO(V)$ & $M_\Gamma$ classically & accidental iso.\\ \hline\hline
0 & 20 & $\SO(2)$ & $\UU(1)$& CM \\ \hline
1 & 19 & $\SO(2,1)$  &$\SL(2,\R)$& Shimura curves\\ \hline
2 & 18 & $\SO(2,2)$  &$\SL(2,\R) \times \SL(2,\R)$ & Hilbert modular surfaces\\ \hline
3 & 17 & $\SO(2,3)$  & $\Sp(2,\R)$ & Siegel $3$-folds\\ \hline 
4 & 16 & $\SO(2,4)$ & $\text{\rm SU}(2,2)$& unitary Shimura $4$-folds\\ \hline
-- & -- & -- & -- & -- \\ \hline
19 & 1 &  \SO(2,19) & -- &moduli of polarized K3's\\ \hline
\end{tabular}
\end{center}

\noindent 
For example, for $n=19$, we have $L' = (2d)$, and we get the moduli space of polarized K3's of polarization degree $2d$. 
Note that $K$ is an even lattice and represents every positive even integer.\hfb
At the other extreme, for $n=0$, we have $\sig(L) =(2,0)$, $\text{rank}\, L'=20$,  and $X$ is a singular K3 surface.\hfb

\subsection{Modular interpretation of the special cycles} 
In the case of families of lattice polarized K3 surfaces $M_\Gamma$ for a lattice $L$ of signature $(2,n)$ as described in the previous 
section,  vectors in $L$ correspond to additional elements of $\Pic(X)$. 
Let $N= \Z^r$, and, for $T\in \Sym_r(\Z)_{>0}$, let $N = N_T$ be the quadratic lattice of signature $(0,r)$ defined by 
$-T$. 

\begin{prop} The codimension $r$ cycle $Z(T)$ can be identified with the locus of objects $(X,\l,j)$
where 
$$j: N_T \hookrightarrow \Pic(X)$$
is a quadratic embedding with
$$j(N_T) \cdot \l(L') = 0.$$
\end{prop}

\noindent Here, if $(X,\l,u)$ is a marked object, then 
$$u\circ j: \Z^r=N_T \hookrightarrow L,\qquad e_i \mapsto x_i,$$
determines an $r$-tuple $\xx\in L^r$ with $T(\xx)=-T$ and the period point of $(X,u)$ lies in $D_{\xx}\subset D(L)$. 

\begin{rem}{\bf(1)} In this construction, we have fixed the basis $\Z^r \simeq N$. 
A change in this basis corresponds to a right multiplication of the row vectors $\xx$ by an element of $\GL_r(\Z)$. \hfb 
{\bf(2)} For $r=1$, we are imposing a single additional class in $\Pic(X)$ and  the $Z(t)$'s are essentially
the Noether-Lefschetz divisors in $M_\Gamma$, cf. \cite{maulik.pand}. 
\end{rem}

\subsection{Some applications:} \hfb

\noindent{\bf I.} Suppose that, for a smooth projective curve $C$, 
$$i_\pi:C\lra M_\Gamma$$
is a morphism corresponding to a family 
$$\pi:X \lra C$$
of $L'$-polarized K3 surfaces. 
Recall that $Z(t) \rightarrow M_\Gamma$ is the locus of collections $(X,\l,j)$, $j\cdot j = -t$,
and consider the fiber product 
$$
\begin{matrix} Z(t)\times_{M_\Gamma}C &\lra & Z(t)\\
\nass
\downarrow&{}&\downarrow\\
\nass
C&\lra&M_\Gamma.
\end{matrix}
$$
Then, for $[C] \in H^{2n-2}_c(M_\Gamma)$,  
\begin{align*}
[Z(t)]\cdot [C] &= \deg \OO_{Z(t)}\vert_{C}.\\
\nass
{}&=\sum_{z\in C} \#\{ (X_z,\l_z,j)\mid j\cdot j = -t, \ j\cdot \l =0\}\qquad\text{(generically)}\\
\nass
{}&=:m(t, X/C) = \text{Noether-Lefschetz number.}
\end{align*}

For example, if the Picard number of a generic member of the family is $20-n$, then $C$ is not contained in any of the 
$Z(t)$'s and the loci in question are all finite 
sets of points. 

\begin{cor}
$$\phi_1(\tau)\cdot [C] = \deg \Cal L\vert_C+\sum_{t\in\Z_{>0}}  m(t,X/C)\, q^t,$$
is an elliptic modular form of weight $\frac{n}2+1$ and level determined by $L$. 
In particular, the numbers $m(t, X/C)$ are the Fourier coefficients of this form. 
\end{cor}

This result is due to Maulik-Pandharipande, \cite{maulik.pand},  where it is derived from Borcherds' Theorem and its significance in Gromov-Witten 
theory is explained. 
\medskip

\noindent{\bf II.}  Similarly, suppose that 
$$\pi:X\lra S, \qquad i_\pi:S \lra M_\Gamma$$
is a family of $L'$-polarized K3 surfaces where $S$ is a projective surface. For $T\in \Sym_2(\Z)_{>0}$, 
$Z(T) \rightarrow M_\Gamma$ is the locus of collections $(X,\l,\bold j)$,  where $\bold j = [j_1,j_2]$ is a pair of classes in $\Pic(X)$ 
orthogonal to $\l(L')$ and with matrix of inner products $\bold j\cdot\bold j = -T$.
Consider   
$$
\begin{matrix} Z(T)\times_{M_\Gamma}S &\lra & Z(T)\\
\nass
\downarrow&{}&\downarrow\\
\nass
S &\lra&M_\Gamma.
\end{matrix}
$$
Then, for the cohomology class $[S] \in H^{2n-4}_c(M_\Gamma)$,  we can define
\begin{align*}
[Z(T)]\cdot [S] &=\sum_{z\in S} \#\{ (X_z,\l_z,\bold j)\mid \bold j\cdot\bold  j = -T, \ \bold j\cdot \l =0\}\qquad\text{(generically)}\\
\nass
{}&=:m(T, X/S) = \text{a higher Noether-Lefschetz number.}
\end{align*}
Here there can be curves in $Z(T)\cap S$ and, in this case, more care must be taken to interpret the intersection number $[Z(T)]\cdot [S]$ . 
By the modularity results above, the $m(T, X/S)$ are Fourier coefficients of a Siegel modular form of genus $2$ and weight $\frac{n}2+1$.
It would be interesting to compute these Siegel modular forms for specific families $X\rightarrow S$,   e.g., ones coming from 
explicit classical geometry. 

\medskip

\noindent{\bf III.}    Suppose that $L$ has signature $(2,2)$ and is anisotropic. 
Then we have a family where $S = M_\Gamma$ is itself a projective surface, so that 
$$\phi_2(\tau,M_\Gamma) \in H^4(M_\Gamma) \overset{\deg}{\isoarrow} \C.$$
 The following result from \cite{kudla.duke} is obtained by combining the results of \cite{KM.IHES} with the extended Siegel-Weil formula of \cite{KR.crelleI}. 
 \begin{theo}\label{thm3.5}  Assume that condition (\ref{assume.3.5}) below holds. Then 
 $$\deg\,\phi_2(\tau,M_\Gamma) = E(\tau,\frac12,L)$$
 where $E(\tau,s,L)$ is a Siegel  Eisenstein series of genus $2$ and weight $2$ associated to $L$, 
evaluated at the Siegel-Weil critical point $s=s_0=\frac12$. 
 \end{theo} 
 
 Note that, for $T>0$, $Z(T)$ is a $0$-cycle, and, when $T \ge 0$ has rank $1$, then $Z(T)^{\text{naive}}$ is a curve on $M_\Gamma$. 
 Thus 
 \begin{multline*}
 \deg\,\phi_2(\tau,M_\Gamma) = \vol(M_\Gamma,\o^2) + \sum_{\text{rank}(T)=1}\vol(Z(T)^{\text{naive}}, \o) \,q^T \\
 \nass
 + \sum_{T>0} \deg(Z(T))\,q^T.
 \end{multline*}
Here, the Siegel-Eisenstein series is defined by
$$E(\tau, s, L)  = \sum_{\gamma\in \Gamma'_\infty\back \Gamma'} \det(c\tau+d)^{\frac{n}2+1} |\det(c\tau+d)|^{s-s_0} \det(v)^{\frac12{s-s_0}}\,\P(\gamma,L),$$
where 
$$\gamma = \begin{pmatrix} a&b \\c&d\end{pmatrix}\in \Gamma' = \Sp_2(\Z).$$
and $\Phi(\gamma,L)$ is a generalized Gauss sum attached to $\gamma$ and $L$, \cite{siegel.book}, \cite{K.Kyoto}.
This series is termwise absolutely convergent for $\Re(s)> \frac32$.
Its value at $s_0=\frac12$ is defined by analytic continuation, \cite{KR.crelleI}.  The main point behind Theorem~\ref{thm3.5} is that, as explained in 
Remark~2.4 (3),  there is a genus $2$ theta function $\theta(\tau,\ph_{KM}^{(2)})$ valued in $A^{(2,2)}(M_\Gamma)$, i.e., in top degree forms, and 
the degree generating series $\deg\,\phi_2(\tau,M_\Gamma)$ is obtained by integrating this form over $M_\Gamma$.  Let $H= \text{\rm O}(V)$ be the orthogonal 
group of $V$ and assume that   there is an open compact subgroup $K\subset H(\A_f)$,  the group of finite ad\`ele points of $H$, such that 
\beq\label{assume.3.5}
H(\A) = H(\Q) H(\R) K, \qquad \Gamma = H(\Q)\cap K.
\eeq
Then the geometric integral of $\theta(\tau,\ph_{KM}^{(2)})$ over 
$$M_\Gamma \simeq H(\Q)\back H(\A)/K_\infty K$$ 
coincides with the integral of a scalar valued theta function on $H(\A)$ over the adelic quotient  $H(\Q)\back H(\A)$  -- cf. \cite{Bints}, section 4,  for a more detailed 
discussion.   The extended Siegel-Weil formula, \cite{KR.crelleI}, \cite{K.Kyoto}, 
identifies the result as a special value at $s=s_0$, perhaps outside the range of absolute convergence,  of a certain Eisenstein series attached to $L$.   

\medskip

\noindent{\bf IV.}  Here is an amusing example.  Let $F=\Q(\sqrt{d})$, $d\in \Z_{>0}$ square free, be a real quadratic field 
with ring of integers $O_F$.  Let $M$ be a projective $O_F$-lattice with a symmetric $O_F$-bilinear form $(\ ,\ )_M$, and 
suppose that the signature of $M$ is given by 
$$\sig(M) = ((2,m), (0,m+2)).$$
Let $L$ be $M$, viewed as a $\Z$-module, with bilinear form $(\ ,\ )_L$ given by 
$$(x,y)_L = \tr_{F/\Q}(x,y)_M.$$
The signature of $L$ is $(2,2m+2)$.  Let $V= L\tt_\Z\Q$ and note that  
$$V\tt_\Q\R = (V\tt_F F)\tt_\Q\R = (V\tt_{F,\s_1}\R)\times (V\tt_{F,\s_2}\R),$$
where $\s_1$ and $\s_2$ are the two real embeddings of $F$. The two factors on the right have signatures $(2,m)$ and $(0,m+2)$ respectively. 
Then there is an embedding
$$D(M) = D(V\tt_{F,\s_1}\R) \hookrightarrow D(V_\R) = D(L).$$
Let $\Gamma_M\subset \Gamma$ be the subgroup of $O_F$-linear isometries of in $\Gamma$. 
Then
$$\Gamma_M\back D(M) \lra \Gamma\back D(L),$$
is an algebraic cycle of codimension $m+2$, and, since the quotient $\Gamma_M\back D(M)$ is compact,  this cycle is projective. 

For example, when $m=2$,  we get a projective surface
$$S \lra M_\Gamma = M_\Gamma(2,6), \qquad [S] \in H^{8}_c(M_\Gamma).$$
It would be nice to have a concrete description of the corresponding family of K3 surfaces
$X \rightarrow S$. 
The Siegel modular generating function for higher Noether-Lefschetz numbers for this family is  
$$\phi_2(\tau, M_\Gamma)\cdot [S] = i^* E_F(\tau,\frac12,\P),\qquad i: \SL_2/\Q \lra \SL_2/F,$$
the pullback of the special value of a genus $2$ Hilbert-Siegel Eisenstein series $E_F(\mathbf\tau,\frac12,\P)$ of weight $(2,2)$ over $F$
to a weight $4$ modular Siegel modular form over $\Q$!  This is again a consequence of the extended Siegel-Weil formula  
together with a seesaw identity, \cite{kudla.seesaw}. Note that, in general, such a pullback will now have lots of 
cuspidal components and so the Fourier coefficients of such a form are essentially more complicated than those of an Eisenstein series
like that occurring in Theorem~\ref{thm3.5}. 

\begin{rem}  A more or less explicit geometric construction of an example for $m=1$ is described in \cite{vG}, section 3.
In this case, $\sig(M) = ((2,1), (0,3))$ and the base of the family will be a Shimura curve $C$ over $F$. 
In this case, the generating series for the Noether-Lefschetz numbers will be an elliptic modular form of weight $3$ 
arising as the pullback for a Hilbert modular 
Eisenstein series of weight $(3/2,3/2)$. 
\end{rem}

\medskip

\noindent{\bf V.}  The previous example can be further generalized.  Suppose that $F$ is a totally real number field with ring of integers $O_F$, 
$|F:\Q|=d>1$, and real embeddings $\s_i:F\hookrightarrow \R$, $1\le i\le d$. 
Let $M$, $( \ , \ )_M$ be a quadratic $O_F$-lattice of rank $m+2$ over $O_F$ and with 
$$\sig(M) = ((2,m),(0,m+2)^{d-1}).$$
% \underset{d-1}{\underbrace{(m+2,0), \dots, (m+2,0)}}).$$
Define a quadratic lattice $L$ as in example {\bf IV}, so that 
$$\sig(L) = (2,d(m+2)-2).$$  
Again setting $V=L\tt_\Z\Q$, we have
$$V\tt_\Q\R =  (V\tt_{F,\s_1}\R)\times (V\tt_{F,\s_2}\R)\times\dots \times (V\tt_{F,\s_d}\R),$$
an embedding
$$D(M) = D(V\tt_{F,\s_1}\R) \hookrightarrow D(V_\R) = D(L),$$
and a projective algebraic cycle 
$$Y=\Gamma_M\back D(M) \lra \Gamma\back D(L) = M_\Gamma,$$
of codimension $(d-1)(m+2)$ and dimension $m$. 
For this construction to fall into the world of K3 moduli, we must require that 
$$2\le d(m+2)\le 21,$$
and hence we have the following little table of the possibilities: \hfb (Note that $N= \dim M_\Gamma = d(m+2)-2$.)

\begin{center}
\begin{tabular}{ | c | c | c | c | c | c | c | c | c | c |}
\hline
$d$ & $2$ & $3$ & $4$ & $5$ & $6$ & $7$ & $8$ & $9$ & $10$\\ \hline
$m=\dim Y$ & $0\le m \le8$ & $0\le m\le 5$ & $0\le m\le 3$ & $0\le m\le 2$ & $0, 1$ & $0, 1$ & $0$ & $0$ & $0$ \\ \hline
$N=\dim M_\Gamma$ & $2\le N \le 18$ & $4\le N\le 19 $ & $6\le N\le 19$ & $8\le N\le 18 $ & $10, 16$ & $12, 19$ & $14$ & $16$ & $18$ \\ \hline
\end{tabular}
\end{center}

It would be interesting to give an account of the families of K3 surfaces over the projective varieties $Y$ occurring here.  

One nice case is the following. 
Let $\kay = \Q(\zeta_{13})$ be the $13$th cyclotomic field and let $F = \Q(\zeta_{13}+\zeta_{13}^{-1})$, so that $|F:\Q|=6$. 
Let $M\tt_{O_F}F$ to be the space of elements of trace $0$ in a quaternion algebra $B/F$ with $B\tt_{F\s_1}\R = M_2(\R)$ and 
$B\tt_{F,\s_i}\R \simeq \mathbb H$, the Hamiltonian quaternions, for $i>1$. 
To specify $B$, we need to choose an additional set $\Sigma_{f}(B)$ of finite places of $F$ with $|\Sigma_f(B)|$ odd. 
The algebra $B$ is determined by the condition that, for a finite place $v$ of $F$, $B\tt_F F_v$ is a division algebra if and only if $v\in \Sigma_f(B)$. 
The simplest choice would be $\Sigma_f(B) = \{ v_{13}\}$, where $v_{13}$ is the unique place above $13$.   We fix a maximal order $O_B$ in $B$ and 
let $M$ be the set of trace $0$ elements in it. 
As quadratic form on $M$,  we take $(x,y) = \tr(xy^\iota)$. 
In this case, $V= L\tt_\Z \Q$ has signature $(2,16)$, and, 
as noted in Remark~\ref{rem3.1} (1), the rational quadratic space $V$ occurs as a summand of $K_\Q$. I have not  
checked\footnote{ If it differs from the primitive lattice $V\cap K$, 
then some `level structure' will have to be introduced.} if the lattice $L$ just defined can occur as $V\cap K$.   In any case, we obtain a Shimura 
curve $C$ embedded in the $16$-dimensional moduli space $M_\Gamma$, and the generating series for the Noether-Lefschetz numbers for the associated 
family of K3 surfaces will be an elliptic modular form of weight $9$ arising as the pullback of a Hilbert modular Eisenstein series of weight 
$(3/2,\dots,3/2)$ for $F$.

\section{Kuga-Satake abelian varieties and special endomorphisms}

In this last section, we give another interpretation of the special cycles $Z(T)$.
For convenience, we change the sign and consider rational quadratic forms of signature $(n,2)$. 

\subsection{The Kuga-Satake construction}

The moduli spaces $M_\Gamma$ of lattice polarized K3 surfaces carry families of abelian varieties 
arising from (a slight variant of) the Kuga-Satake construction, which we now review.  

For a quadratic lattice $L$ of signature $(n,2)$ and for $V=L_\Q$, let
$$C=C(V)= C^+(V)\oplus C^-(V)$$
be the Clifford algebra of $V$. 
For a basis $v_1,\dots, v_{n+2}$ of $V$, let
$\delta = v_1\cdots v_{n+2}\in C(V)$, and let $\kay = \Q(\delta)$ be the discriminant field.
For $n$ even, $\delta\in C^+(V)$,  $\kay$ is the center of $C^+(V)$, and $\delta$ anticommutes with elements of $C^-(V)$. 
For $n$ is odd,  $\delta\in C^-(V)$ and $\kay$ is central in all of $C(V)$. 
Now
$$
\OC=C(L)=\text{Clifford algebra of $L$}
$$ 
is an order in $C$ and $\OK = \kay \cap C(L)$ is an order in $\kay$.
We obtain a real torus
$$A^{\top}=C(L_\R)/C(L)$$
of dimension $2^{n+2}$ with an action
$$\iota: C(L) \tt_\Z \OK  \lra \End(A^\top), \qquad \iota(c\tt \a): x \longmapsto c x \a,$$
via left-right multiplication. Let
$$
G = \GSpin(V) =\{ g\in C^+(V)^\times \mid g V g^{-1} = V\}
$$
and, for an element 
$
a\in C^+(L)\cap C(V)^\times$ with $a^\iota = -a$,
define an alternating form on $C$ by
$$
\gs{x}{y}=\gs{x}{y}_a = \tr(a x y^\iota).
$$
Here $x\mapsto x^\iota$ is the main involution on $C(V)$, determined by the conditions 
$(xy)^\iota= y^\iota x^\iota$ for all $x$, $y\in C$, and $v^\iota = v$ for all $v\in V \subset C^-(V)$. 
Note that for $g\in G$, 
$$\gs{xg}{yg} = \nu(g)\,\gs{x}{y}, \qquad \nu(g) = g g^\iota = \text{spinor norm.}
$$

A variation of complex structures on $A^\top$ is defined as follows. 
For an oriented negative $2$-plane $z\in D(L)$ in $V_\R$, let $z_1$, $z_2$ be a properly oriented 
orthonormal basis and let $j_z = z_1z_2\in C^+(V_\R)$. This element is independent of the choice of basis. Note that $j_z^2=-1$, so that 
right multiplication by $j_z$ defines a complex structure on $C(V_\R)$,  and hence we have a complex torus
$A_z = (A^{\top}, j_z)$ for each $z\in D(L)$.  Since the action of $C(L)\tt \OK$ 
commutes with the right multiplication by $j_z$, we have 
\beq\label{endos}
\iota:\OC \tt_\Z\OK \lra \End(A_z)
\eeq
with
$$\iota(c\tt \a)^* = \iota(c^*\tt \a^\iota),\qquad c^* = a c^\iota a^{-1},$$
for the Rosati involution determined by $\gs{\ }{\ }$. 
If $\gamma\in \Gamma \subset \Gamma_L$, and\footnote{In general, there is an orbifold issue here that can be eliminated by introducing a suitable 
level structure.} if 
$\tilde{\gamma} \in O_C^\times$ is an element mapping to $\gamma$ under the natural homomorphism 
$\text{\rm GSpin(V)} \rightarrow \SO(V)$, then right multiplication by $\tilde\gamma$ on $C(L_\R)$ 
induces an isomorphism of the complex tori  
$A_{\gamma(z)}$ and $A_z$, equivariant for the action of $O_C\tt\OK$. 
Finally,  fix a rational splitting $V = V^+ + V^-$ of signature $(n,0) + (0,2)$ and let $a_1$, $a_2$ be a $\Z$-basis 
for the negative definite lattice $L^- = L \cap V^-$. Let $a = a_1a_2$.  With this choice of $a$, the form 
$\gs{\ }{\ }_a$ is $\Z$-valued on $C(L)$ and defines a Riemann form on each complex torus $A_z$, i.e., the form $\gs{x j_z}{y}_a$ is symmetric and definite
(hence positive definite on one of the connected components of $D(L)$). 
Thus, $M_\Gamma$ carries a family of polarized abelian varieties:
$$\pi:\text{KS}(L) \lra M_\Gamma,\qquad (A_z, \iota, \l), \qquad \dim A_z = 2^{n+1}$$
with endomorphisms (\ref{endos}).  Note that these abelian varieties are $\Z_2$-graded, 
$$A_z = A^+_z\oplus A^-_z,$$
since the construction
respects the decomposition $C(L) = C^+(L)\oplus C^-(L)$. 
\begin{rem}
(1)  For $n\le 4$, $M_\Gamma$ is the moduli space for such PE type abelian varieties and is a Shimura variety of PEL type\footnote{We could add a level structure here.}. \hfb
(2) For general $n \ge 5$, the abelian varieties in the family are not characterized by the given PE; more Hodge classes are required, 
and $M_\Gamma$ is a Shimura variety of Hodge type. 
The family $\pi$ corresponds to a morphism
$$i_\pi: M_\Gamma \lra A_g, \qquad g= 2^{n+1}.$$
\end{rem}

\subsection{Special cycles and special endomorphisms}\hfb

Special endomorphisms arise as follows.  For $x\in L$, 
let 
$$r_x:A^\top \lra A^\top = C(L_\R)/C(L),$$
be the endomorphism induced by right multiplication by $x$.  Note that this endomorphism 
has degree $1$ with respect to the $\Z_2$-grading, and that, since $x^2 = Q(x)$ in $C(L)$, 
$$r_x^2= [Q(x)],$$
is just multiplication by the integer $Q(x) = (x,x)\in \Z$. 
Moreover, 
\beq\label{special1}
r_x\circ \iota(c\tt b) = 
\begin{cases} \iota(c\tt b)\circ r_x &\text{for $n$ odd}\\
\nass
\iota(c\tt b^\s)\circ r_x &\text{for $n$ even.}
\end{cases}
\eeq
and the adjoint of $r_x$ with respect to $\gs{\ }{\ }$ is
\beq\label{special2}
(r_x)^*= r_{x^\iota} = r_x.
\eeq
\begin{defn}
For a given $z\in D(L)$, the endomorphism $r_x$ is said to be a special endomorphism of the abelian variety $A_z$ when it is holomorphic, 
i.e., when $r_x\in \End(A_z)$. 
\end{defn}

\begin{lem} Let $L(A_z)$ be the space of special endomorphisms of $A_z$, with quadratic form defined by 
$j^2 = Q(j)$. Then 
$$L(A_z) = L\cap z^\perp$$
is an isometry.
\end{lem}
Note that, $r_x$ is a special endomorphism of $A_z$ precisely when $x$ commutes with $j_z$. 
Since conjugation by $j_z$ induces the endomorphism of $V$ that is $-1$ on $z$ and $+1$ on $z^\perp$, it follows that
$r_x$ is a special endomorphism precisely when $x\in z^\perp$, i.e., $z\in D_x$. 

\begin{cor} For $T\in \Sym_r(\Z)_{>0}$, the special cycle $Z(T)$ is the image in $M_\Gamma$ of
$$\tilde Z(T) = \{ (A_z, \bold j)\mid \bold j\in L(A_z)^r, \  (\bold j, \bold j) =T\ \}.$$
\end{cor}

In fact, to give the right definition of special endomorphism, one should start with an object 
$(A,\iota, \l, \e, \xi)$ where $A$ is an abelian variety (scheme) with $O_C\tt\OK$-action $\iota$, a $\Z/2\Z$-grading $\e$ and 
additional Hodge tensors $\xi$. A special endomorphism is then an element $j\in \End(A)$ of degree $1$ with respect to the 
$\Z/2\Z$-grading satisfying conditions (\ref{special1}) and (\ref{special2}), with an additional compatibility with respect to 
the Hodge tensors $\xi$. 
For $0\le n\le 3$, where the classes $\xi$ can be expressed in terms of endomorphisms, such a definition of special 
endomorphism was given in \cite{kry.tiny} ($n=0$), \cite{KRY.book} ($n=1$), \cite{krcrelle} ($n=2$), and \cite{krsiegel} ($n=3$) 
and used to define special cycles in integral models of the associated PEL type Shimura varieties.  The results of these papers
concerning the intersection number of integral special cycles and Fourier coefficients of modular forms should have some
consequences in the arithmetic theory of K3 surfaces. To reveal it one would have to relate these integral models to the integral models of 
moduli spaces of K3 surfaces developed by Rizov, \cite{rizov.models}, \cite{rizov.KS}. 
The extension to the case of general $n$, i.e., to integral models of Shimura varieties attached to $\text{GSpin}(n,2)$
is the subject of ongoing work of Andreatta-Goren, \cite{AG.oberwolfach} and of Howard-Madapusi, \cite{HM}.

\begin{rem}
In the references that follow, I have made no 
attempt at completeness and apologize in advance for the many omitted citations. 
\end{rem}

\end{document}